\documentclass{amsart}
\usepackage{amsmath, amsthm, amscd, amsfonts, amssymb, graphicx, color}
\usepackage[ plainpages]{hyperref}

\newtheorem{theorem}{Theorem}[section]
\newtheorem{lemma}[theorem]{Lemma}

\theoremstyle{definition}

\theoremstyle{remark}
\newtheorem{remark}[theorem]{\textbf{Remark}}

\usepackage[all]{xy}
\numberwithin{equation}{section}

\begin{document}
\setcounter{page}{1}

\title[ Toeplitz Operators on the Harmonic Bergman Space ]
{ On the Commuting Problem of Toeplitz Operators on the Harmonic Bergman Space }

\author[ H. Iqtish, I. Louhichi, and A. Yousef]{ Hasan Iqtaish$^{1}$, Issam Louhichi$^{2}$ and Abdelrahman Yousef$^{3}$ }

\address{$^{1}$ Department of Mathematics and Statistics, College of Arts and Sciences, American University of Sharjah, Sharjah, UAE.} \email{\textcolor[rgb]{0.00,0.00,0.84}{b00101964@aus.edu}}

\address{$^{2}$ Department of Mathematics and Statistics, College of Arts and Sciences, American University of Sharjah, Sharjah, UAE.} \email{\textcolor[rgb]{0.00,0.00,0.84}{ilouhichi@aus.edu}}

\address{$^{3}$ Department of Mathematics and Statistics, College of Arts and Sciences, American University of Sharjah, Sharjah, UAE.} \email{\textcolor[rgb]{0.00,0.00,0.84}{afyousef@aus.edu}}

\begin{abstract} In this paper, we provide a complete characterization of bounded Toeplitz operators $T_f$ on the harmonic Bergman space of the unit disk, where the symbol $f$ has a  polar decomposition truncated above, that commute with $T_{z+\bar{g}}$, for a bounded analytic function $g$.

  \end{abstract}\maketitle

\maketitle
\section{\textbf{Introduction}}
Let \(dA = r \, dr \, \frac{d\theta}{\pi}\), where \((r, \theta)\) are the polar coordinates, denote the normalized Lebesgue area measure on the unit disk \( \mathbb{D} \),  so that \( \mathbb{D} \) has measure 1. The space \( L^2(\mathbb{D}, dA) \) consists of all Lebesgue square-integrable functions on \( \mathbb{D} \) and forms a Hilbert space with the inner product
\[
\langle f, g \rangle = \int_{\mathbb{D}} f(z) \overline{g(z)} \, dA(z).
\]

The \textbf{harmonic Bergman space}, denoted by \( L^2_h(\mathbb{D}) \), is the closed subspace of \( L^2(\mathbb{D}, dA) \) comprising all complex-valued \( L^2 \)-harmonic functions on \( \mathbb{D} \). Let \( Q \) represent the orthogonal projection from \( L^2(\mathbb{D}, dA) \) onto \( L^2_h(\mathbb{D}) \). This projection is given by the integral operator
\[
Q f(z) = \int_{\mathbb{D}} \left( \frac{1}{(1 - z\bar{w})^2} + \frac{1}{(1 - z\bar{w})^2} - 1 \right) f(w) \, dA(w), \quad z \in \mathbb{D},
\]
for any \( f \in L^2(\mathbb{D}, dA) \). It is well-known that $Q$ is bounded from $L^2(\mathbb{D},dA)$ onto $L^2_h(\mathbb{D})$. 

For a function \( u \in L^1(\mathbb{D}, dA) \), we define the \textbf{Toeplitz operator} \( T_u \) with symbol \( u \) on \( L^2_h \) by
\begin{equation}
T_u f = Q(u f)
\end{equation}
for $f \in L^2_h(\mathbb{D})$, provided that the product $uf$ is in $L^2(\mathbb{D}, dA)$. This operator is densely defined on the polynomials and not bounded in general. However, if $u$ is bounded on $\mathbb{D}$, then $T_u$ is bounded and $||T_u||\leq ||u||_{\infty}$.

A symbol \(u\) is called quasihomogeneous of degree \(p\), where \(p\) is an integer, if it can be expressed in the form \(u(re^{i\theta}) = e^{ip\theta}\phi(r)\), where  \(\phi\) is a radial function. In this case, the associated Toeplitz operator \(T_u\)  is called a quasihomogeneous Toeplitz operator of degree \(p\). The study of these operators is motivated by the structural decomposition of \(L^2(\mathbb{D}, dA)\), which can be written as \(L^2(\mathbb{D}, dA) = \bigoplus_{k \in \mathbb{Z}} e^{ik\theta} \mathcal{R}\), where \(\mathcal{R}\) denotes the space of square-integrable radial functions on \([0,1)\) with respect to the measure \(rdr\).  This decomposition implies that any function \(f \in L^2(\mathbb{D}, dA)\) admits a polar decomposition \(f(z) = f(re^{i\theta}) = \sum_{k \in \mathbb{Z}} e^{ik\theta} f_k(r)\), where each \(f_k(r)\) is a radial function. Moreover, we say that $f$ is truncated above if its polar decomposition is of the form $f(re^{i\theta})=\sum_{k=-\infty}^Ne^{ik\theta}f_k(r)$, for some integer $N$.

Our focus is on identifying the conditions that characterize the symbols of commuting Toeplitz operators on $L^2_h(\mathbb{D})$. This problem has been extensively explored in the contexts of the classical Hardy space and the analytic Bergman space over the years. The study of Toeplitz operators on $L^2_h(\mathbb{D})$ exhibits significant differences compared to their counterparts on the analytic Bergman space and remains less understood. However, there has been growing
interest in investigating this issue within the framework of the harmonic Bergman space. For instance, Choe and Lee \cite{choe} established that two analytic Toeplitz operators on $L^2_h$  i.e., Toeplitz operators with analytic symbols,  commute if and only if their symbols, along with the constant function $1$, are linearly dependent. Subsequent works such as \cite{choe1} and \cite{ding} demonstrated that an analytic Toeplitz operator and a co-analytic Toeplitz operator on $L^2_h(\mathbb{D})$ can commute only if at least one of their symbols is a constant function. In \cite{dong2}, the conditions under which the product of two quasihomogeneous Toeplitz operators remains a Toeplitz operator were investigated. Building upon this, the work in \cite{dong3} delved into the commuting problem for quasihomogeneous Toeplitz operators on $L^2_h(\mathbb{D})$, where the authors characterized the commuting Toeplitz operators with quasihomogeneous symbols. In addition, they showed that a Toeplitz operator with an analytic or co-analytic monomial symbol commutes with another Toeplitz operator only in the trivial case. For further results on commuting Toeplitz operators in harmonic Bergman spaces, the reader may consult \cite{chen, choe, choe1, choe2, ding, ding2, dong1, dong2, dong3, dong4, guan, guo, kong, le, lz, lzf} and the references therein.\\

The primary goal of our study is to characterize a special class of commuting Toeplitz operators acting on $L^2_h(\mathbb{D})$. More specifically, we characterize all Toeplitz operators with truncated above symbols that commute with the Toeplitz operator $T_u$, whose symbol is the harmonic function  $u(z)=z+\overline{g(z)},\   \textrm{ where } g(z)=\sum_{n=0}^{\infty}a_nz^n \textrm{ is a bounded analytic function on } \mathbb{D}$.  

One of the main challenges in this problem arises from the interplay between multiplication operators induced by the symbols and the projection onto the harmonic Bergman space. Unlike the analytic Bergman space, where the Bergman projection has an explicit integral representation, the harmonic Bergman projection introduces additional complexities that make computing Toeplitz operator products more difficult. Consequently, many classical results from the analytic setting do not directly extend to the harmonic case, necessitating the development of new techniques and approaches.

To structure our analysis effectively, this paper is organized as follows. Section 2 presents key preliminary results essential for proving the main theorem. Section 3  formally states the main result. Finally, Section 4 is devoted to its proof, which is divided  into several lemmas to enhance clarity and systematically manage the technical details.
\section{\textbf{Tools}}

The \textbf{Mellin transform} \( \widehat{\phi} \) of a radial function \( \phi \in L^1([0,1), r \, dr) \) is given by
\[
\widehat{\phi}(z) = \int_0^1 \phi(r) r^{z-1} \, dr.
\]
It is well-known that for such functions, the Mellin transform is bounded in the right half-plane \( \{z \in \mathbb{C} : \Re z \geq 2\} \) and is analytic in  \( \{z \in \mathbb{C} : \Re z > 2\} \).

The following lemma describes the action of quasihomogeneous Toeplitz operators on elements of the orthogonal basis of \( L^2_h(\mathbb{D}) \). See   \cite[Lemma 2.1, p.~1767]{dong2}.

\begin{lemma}\label{mellin}Let \( k \in \mathbb{Z} \) and let \( \phi \) be a radial in $L^{1}([0,1), rdr)$. Then, for each \( n \in \mathbb{N} \), the Toeplitz operator $T_{e^{ik\theta}\phi}$ satisfies
\[
T_{e^{ik\theta} \phi}(z^n) =
\begin{cases}
2(n + k + 1) \widehat{\phi}(2n + k + 2) z^{n+k}, & \text{if } n \geq -k, \\
2(-n - k + 1) \widehat{\phi}(-k + 2) \overline{z}^{-n-k}, & \text{if } n < -k.
\end{cases}
\]

Similarly,
\[
T_{e^{ik\theta} \phi}(\bar{z}^n) =
\begin{cases}
2(n - k + 1) \widehat{\phi}(2n - k + 2)\overline{z}^{n-k}, & \text{if } n \geq k, \\
2(k - n + 1) \widehat{\phi}(k + 2) z^{k-n}, & \text{if } n < k.
\end{cases}
\]
\end{lemma}

A fundamental result states that the Mellin transform of a function is uniquely determined by its values on an arithmetic sequence of integers. This is formalized in the following classical theorem \cite[p.~102]{Rem}.

\begin{theorem}\label{rem}
Let \( f \) be a bounded analytic function in the right half-plane \( \{z \in \mathbb{C} : \Re z > 0\} \) that vanishes at an infinite sequence of distinct points \( d_1, d_2, \dots \) satisfying:
\begin{itemize}
\item[(i)] \( \inf \{ |d_n| \} > 0 \), and
\item[(ii)] \( \sum\limits_{n \geq 1} \Re\left( \frac{1}{d_n} \right) = \infty \).
\end{itemize}
Then \( f \) must be identically zero on \( \{z \in \mathbb{C} : \Re z > 0\} \).
\end{theorem}

Another important result we frequently use is the following  lemma. See \cite[Lemma 7, p.~1727]{le2}.

\begin{lemma}\label{nev}
If a meromorphic function in a right half-plane belongs to the Nevanlinna class and is periodic, then it must be constant.
\end{lemma}

The following lemma is crucial for the proof of the main result and can be deduced from  \cite[Theorem 3.8, p.~1278]{dong3}.
\begin{lemma}
\label{anal}
Let $f(re^{i\theta})=e^{ip\theta}\phi(r)$ be a quasihomogenous symbol, where $p\in \mathbb{Z}_+$ and $\phi(r)\in L^1([0,1),rdr)$. If $T_fT_{z^n}=T_{z^n}T_f$ for $n\geq 1$ integer, then $\phi(r)=Cr^p$. In other words, $f$ must be analytic of the form $f(z)=Cz^p$.  
\end{lemma}
\begin{remark}
\label{periodic}
The following observations will be useful in our proves:
\begin{itemize}
\item[1)] A straightforward calculation shows that 
\[
\widehat{r^n}(z) = \frac{1}{z+n}, \quad \text{for } n \in \mathbb{Z},
\] and
$$ \widehat{r^a\ln (r)^b}=\dfrac{(-1)^b b!}{(a+z)^b}, \textrm{where } a>0\textrm{ and } b \textrm{ is a nonnegative integer.}$$
\item[2)] Regarding Theorem \ref{rem}, we apply it in the following setting: Suppose \( (n_k)_k \) is an arithmetic sequence of positive integers and that, for some radial function \( \phi \), we have \( \widehat{\phi}(n_k) = 0 \) for all \( k \). By Theorem \ref{rem}, this forces \( \widehat{\phi} \) to be identically zero in the right half-plane, implying that \( \phi \) itself must vanish there as well.
\item[3)] Lemma \ref{nev} is a key tool in our arguments. In  our proofs, we encounter functional equations of the form
\[
F(z + p) - F(z) = G(z + p) - G(z),
\]
where \( \Re(z) > 0 \), \( p \) is an integer, and \( F \) and \( G \) are bounded analytic functions in the right half-plane. Applying Lemma \ref{nev}, we conclude that \( F(z) = C + G(z) \) for some constant \( C \).
\end{itemize}
\end{remark}
\section{Main Result}
Given a symbol $u(z)=z+\overline{g(z)}$, where  $g(z)=\sum_{n=1}^{\infty}a_nz^n$ is a bounded analytic function on  $\mathbb{D},$ we aim to characterize all symbols of the form (i.e., symbols whose polar decomposition is truncated above)
	$$f(re^{i\theta})=\sum_{n=-\infty}^Ne^{in\theta}f_n(r),\ N\geq 1,$$ in $L^1(\mathbb{D},dA)$ for which the associated Toeplitz operators $T_f$ are bounded and commute with $T_u$.  It is understood here that $f_N\neq 0$.
	We recall that $T_f$ commutes with $T_u$ if and only if  
	\begin{equation}\label{commute} T_fT_u(z^k)=T_uT_f(z^k)\end{equation} and \begin{equation}\label{commute1} T_fT_u(\bar{z}^k)=T_uT_f(\bar{z}^k)\end{equation} for all vectors $z^k$ and $\bar{z}^k$ in the orthogonal basis of $L^2_h(\mathbb{D})$. 
	
	Our main theorem can be stated as follows.
	\begin{theorem}
\label{main}
		Let $u(z)=z+\sum_{l=1}^{\infty}\bar{a_l}\bar{z}^l$. If there exists a nonzero function $f$ of the form $f(re^{i\theta})=\sum_{n=-\infty}^Ne^{in\theta}f_n(r)$,  with $N\geq 1,$ such that $T_f$ commutes with $T_u$, then $T_f$ is a polynomial of degree at most one in $T_u$. In other words, there exist constants $C_1, C_0$ such that $T_f=C_1T_u+C_0I$, where $I$ denotes the identity operator.
	\end{theorem}
\section{Proof of the main result}	
The proof of our main result is quite lengthy and involves intricate computations. To enhance clarity and readability, we have structured the proof into several lemmas. The first lemma establishes that the highest degree $N$ of $f$ in Theorem\ref{main} cannot exceed $3$. However, we will later demonstrate that $N$ must, in fact, be equal to $1$.
\begin{lemma}
\label{degree}
Under the hypothesis Theorem \ref{main}, we have that  $N\leq 3$.
\end{lemma}

\begin{proof}
The term in $z$ of degree $n+N+1$ in equation \eqref{commute} appears on both sides, originating from
$T_{e^{iN\theta}f_N}T_{z}(z^n)$ and $T_{z}T_{e^{iN\theta}f_N}(z^n)$, respectively. Therefore, we have $$T_{e^{iN\theta}f_N}T_{z}(z^n)=T_{z}T_{e^{iN\theta}f_N}(z^n)$$ for every $n$. By Lemma\ref{anal}, this implies that  $e^{iN\theta}f_N=C_Nz^N$ for some constant $C_N$. Similarly, the term in $z$ of degree $n+N$ appears on both sides only from  $T_{e^{i(N-1)\theta}f_{N-1}}T_{z}(z^n)$ and $T_{z}T_{e^{i(N-1)\theta}f_{N-1}}(z^n)$. Applying Lemma \ref{anal} again, we conclude that $e^{i(N-1)\theta}f_{N-1}$ is analytic and satisfies $e^{i(N-1)\theta}f_{N-1}=C_{N-1}z^{N-1}$ for some constant $C_{N-1}$.

Next, the term in $z$ of degree $n+N-1$ comes from
\begin{eqnarray*}
\left(T_{C_Nz^N}T_{\bar{a_1}\bar{z}}+T_{e^{i(N-2)\theta}f_{N-2}}T_{z}\right)(z^n)&=&\left(T_{\bar{a_1}\bar{z}}T_{C_Nz^N}+T_{z}T_{e^{i(N-2)\theta}f_{N-2}}\right)(z^n)
\end{eqnarray*}
Using Lemma \ref{mellin}, the previous equation implies
\begin{eqnarray*}
C_N\bar{a_1}\frac{2n}{2n+2}+(2n+2N)\widehat{f}_{N-2}(2n+N+2)&=&C_N\bar{a_1}\frac{2n+2N}{2n+2N+2}\\&+&(2n+2N-2)\widehat{f}_{N-2}(2n+N),
\end{eqnarray*}
which is equivalent to
\begin{eqnarray*}
2(n+N)\widehat{f}_{N-2}(2n+N+2)-2(n+N-1)\widehat{f}_{N-2}(2n+N)&=&C_N\bar{a_1}\frac{2n+2N}{2n+2N+2}\\
&-&C_N\bar{a_1}\frac{2n}{2n+2}.
\end{eqnarray*}
We complexify the above equation by letting $z=2n$ and we use Remark \ref{periodic} to obtain
\begin{eqnarray}
\label{periodic1}
(z+2N)\widehat{f}_{N-2}(z+N+2)-(z+2N-2)\widehat{f}_{N-2}(z+N)&=&C_N\bar{a_1}\frac{z+2N}{z+2N+2}\nonumber\\ &-&C_N\bar{a_1}\frac{z}{z+2}.
\end{eqnarray}
Define $F(z)=(z+2N-2)\hat{f}_{N-2}(z+N)$ and $\displaystyle G(z)=C_N\bar{a_1}\sum_{i=0}^{N-1}\frac{z+2i}{z+2i+2}$. Then equation \eqref{periodic1} becomes
$$F(z+2)-F(z)=G(z+2)-G(z).$$
Thus, Remark \ref{periodic} implies the existence of a constant $C_{N-2}$ such that $F(z)-G(z)=C_{N-2}$. Equivalently,
\begin{equation}\label{N-2}
(z+2N-2)\hat{f}_{N-2}(z+N)=C_{N-2}+C_N\bar{a}_1\sum_{i=0}^{N-1}\frac{z+2i}{z+2i+2}
\end{equation}
or
\begin{eqnarray*}
\widehat{r^Nf_{N-2}}(z)&=&\frac{C_{N-2}}{z+2N-2}+C_N\bar{a}_1\sum_{i=0}^{N-1}\frac{z+2i}{(z+2N-2)(z+2i+2)}\\
&=& \frac{C_{N-2}}{z+2N-2}+C_N\bar{a}_1\left[\frac{1}{z+2N}+\sum_{i=0}^{N-2}\frac{z+2i}{(z+2N-2)(z+2i+2)}\right]\\
&=& \frac{C_{N-2}}{z+2N-2}+C_N\bar{a}_1\Bigg[\frac{1}{z+2N}+\frac{z+2N-4}{(z+2N-2)^2}\\
&+&\sum_{i=0}^{N-3}\frac{z+2i}{(z+2N-2)(z+2i+2)}\Bigg]\\
&=& \frac{C_{N-2}}{z+2N-2}+C_N\bar{a}_1\Bigg[\frac{1}{z+2N}+\frac{1}{z+2N-2}-\frac{2}{(z+2N-2)^2}\\
&+&\sum_{i=0}^{N-3}\left(\frac{2-2N+2i}{(4-2N+2i)(z+2N-2)}-\frac{2}{(2N-2i-4)(z+2i+2)}\right)\Bigg]\\
&=&  C_{N-2}\widehat{r^{2N-2}}(z)+C_N\bar{a}_1\Bigg[\widehat{r^{2N}}(z)+\widehat{r^{2N-2}}(z)+\widehat{2r^{2N-2}\ln r}(z)\\
&+&\sum_{i=0}^{N-3}\left(\frac{2-2N+2i}{4-2N+2i}\widehat{r^{2N-2}}(z)-\frac{2}{2N-2i-4}\widehat{r^{2i+2}}(z)\right)\Bigg]
\end{eqnarray*}
Therefore, Remark \ref{periodic} implies
\begin{eqnarray*}f_{N-2}(r)&=& C_{N-2}r^{N-2}+C_N\bar{a}_1\Bigg[r^{N}+r^{N-2}+2r^{N-2}\ln r+ \sum_{i=0}^{N-3}\Bigg(\frac{2-2N+2i}{4-2N+2i}r^{N-2}\\
	&-&\frac{2}{2N-2i-4}r^{2i+2-N}\Bigg)\Bigg]\end{eqnarray*}
Observe that $f_{N-2}$ belongs $L^1([0,1),rdr)$ if and only if  $2i+2-N+1\geq 0$, which simplifies to $N\leq 2i+3$ for all $i=0,1,2,\ldots,N-3$. Consequently, this condition must hold for $i=0$, which leads to  $N\leq 3$.
\end{proof}
\begin{remark}
\label{rmk1}
Lemma \ref{degree} implies the following:
\begin{enumerate}
\item $f_3(r)=C_3r^3$.\\
\item $f_2(r)=C_2r^2$.\\
\item To find $f_1(r)$, we plug $N=3$ in equation \eqref{N-2} to obtain
\begin{eqnarray*}
\widehat{r^3f_{1}}(z)&=&\frac{C_{1}}{z+4}+C_3\bar{a}_1\sum_{i=0}^{2}\frac{z+2i}{(z+4)(z+2i+2)}\\
                     &=& \frac{C_{1}}{z+4}+C_3\bar{a}_1\left[\frac{z}{(z+4)(z+2)}+\frac{z+2}{(z+4)^2}+\frac{z+4}{(z+6)(z+4)}\right]\\
&=& \frac{C_{1}}{z+4}+C_3\bar{a}_1\left[\frac{-1}{z+2}+\frac{3}{z+4}+\frac{-2}{(z+4)^2}+\frac{1}{z+6}\right]\\
&=&C_{1}\widehat{r^{4}}(z)+C_3\bar{a}_1\left[\widehat{r^{6}}(z)+3\widehat{r^{4}}(z)+2\widehat{r^{4}\ln r}(z)-\widehat{r^{2}}(z)\right].
\end{eqnarray*}
Hence, Remark \ref{periodic} yields $$f_1(r)=C_1r+C_3\bar{a}_1\left[r^3+3r+2r\ln r-\frac{1}{r}\right].$$
\end{enumerate}
\end{remark}
So far, using Lemma \ref{degree} and Remark \ref{rmk1}, we have established that any Toeplitz operator with symbol $f(re^{i\theta})=\displaystyle{\sum_{k=-\infty}^{N}e^{ik\theta}f_{k}(r)}$ that commutes with $T_{z+\overline{g}}$, where $g(z)=\displaystyle{\sum_{l=1}^{\infty}a_lz^l}$ is a bounded analytic function on $\mathbb{D}$, must take the form
$$f(re^{i\theta})=C_3z^3+C_2z^2+e^{i\theta}\left(C_1r+C_3\bar{a}_1\left[r^3+3r+2r\ln r-\frac{1}{r}\right]\right)+\displaystyle{\sum_{k=-\infty}^{0}e^{ik\theta}f_{k}(r)}.$$

In the following lemmas, we compute the exact expressions of $f_0(r)$, $f_{-1}(r)$, and $f_{-2}(r)$.

\begin{lemma}
\label{f0}
Under the hypothesis of Theorem \ref{main}, we have that $$f_0(r)=C_0+C_2\bar{a}_1\left[1+2\ln r+r^2\right]+C_3\bar{a}_2\left[4\ln r+2r^2+r^4\right].$$
\end{lemma}
\begin{proof}
For $n\geq 1$, the term $z^{n+1}$ in $T_{f}T_{z+\overline{g}}(z^n)=T_{z+\overline{g}}T_{f}(z^n)$ appears both sides only from the expressions $$\left(T_{f_0}T_z+T_{e^{2i\theta}f_2}T_{\bar{a}_1\bar{z}}+T_{e^{3i\theta}f_3}T_{\bar{a}_2\bar{z}^2}\right)(z^n)$$ and $$\left(T_zT_{f_0}+T_{\bar{a}_1\bar{z}}T_{e^{2i\theta}f_2}+T_{\bar{a}_2\bar{z}^2}T_{e^{3i\theta}f_3}\right)(z^n),$$ respectively. 
Thus, applying Lemma \ref{mellin}, we obtain
\begin{eqnarray*}
(2n+4)\widehat{f}_0(2n+4)+C_2\bar{a}_1\frac{2n}{2n+2}+C_3\bar{a}_2\frac{2n-2}{2n+2}&=&(2n+2)\widehat{f}_0(2n+2)+C_2\bar{a}_1\frac{2n+4}{2n+6}\\
&+&C_3\bar{a}_2\frac{2n+4}{2n+8},
\end{eqnarray*}
which can be written as
$$(2n+4)\widehat{f}_0(2n+4)-(2n+2)\widehat{f}_0(2n+2)=C_2\bar{a}_1\left[\frac{2n+4}{2n+6}-\frac{2n}{2n+2}\right]+C_3\bar{a}_2\left[\frac{2n+4}{2n+8}-\frac{2n-2}{2n+2}\right].$$
We complexify the above equation by considering $z=2n$ and we use Remark \ref{periodic} to obtain
\begin{eqnarray*}
\label{complex1}
(z+6)\widehat{f}_0(z+6)-(z+4)\widehat{f}_0(z+4)&=&C_2\bar{a}_1\left[\frac{z+6}{z+8}-\frac{z+2}{z+4}\right]+C_3\bar{a}_2\left[\frac{z+6}{z+10}-\frac{z}{z+4}\right].
\end{eqnarray*}
This equation can be expressed in the form $$F(z+2)-F(z)=G(z+2)-G(z),$$ where $F(z)=(z+4)\widehat{f}_0(z+4)$ and $\displaystyle G(z)=C_2\bar{a}_1\sum_{i=0}^1\frac{z+2i+2}{z+2i+4}+C_3\bar{a}_2\sum_{i=0}^2\frac{z+2i}{z+2i+4}$. By Remark \ref{periodic},  it follows that there exists a constant $C_0$ such that $F(z)-G(z)=C_0$. Hence, we have
\begin{eqnarray*}
\widehat{r^4f_0}(z)&=&\frac{C_0}{z+4}+C_2\bar{a}_1\Bigg[\frac{z+2}{(z+4)^2}+\frac{1}{z+6}\Bigg]+C_3\bar{a}_2\Bigg[\frac{z}{(z+2)^2}+\frac{z+2}{(z+4)(z+6)}\\
&+&\frac{1}{z+8} \Bigg]\\
                   &=& \frac{C_0}{z+4}+C_2\bar{a}_1\Bigg[\frac{1}{z+4}-\frac{2}{(z+4)^2}+\frac{1}{z+6}\Bigg]+C_3\bar{a}_2\Bigg[\frac{1}{z+4}-\frac{4}{(z+4)^2}\\&+&\frac{2}{z+6}-\frac{1}{z+4}+\frac{1}{z+8} \Bigg]\\
                   &=& C_0\widehat{r^4}(z)+C_2\bar{a}_1\left[\widehat{r^4}(z)+2\widehat{r^4\ln r}(z)+\widehat{r^6}(z)\right]+C_3\bar{a}_2\Big[4\widehat{r^4\ln r}(z)\\&+&2\widehat{r^6}(z)+\widehat{r^8}(z)\Big].
\end{eqnarray*}
Therefore, $$f_0(r)=C_0+C_2\bar{a}_1\left[1+2\ln r+r^2\right]+C_3\bar{a}_2\left[4\ln r+2r^2+r^4\right].$$
\end{proof}

Next, we proceed to compute  the radial function $f_{-1}$.
\begin{lemma} 
\label{f-1} Under the hypothesis of Theorem \ref{main}, we have that
 \begin{eqnarray*}f_{-1}(r)&=&\frac{C_{-1}}{r}+C_1\bar{a}_1r+C_3\bar{a}_1^2\left[3r+2r\ln r+r^3\right]+C_2\bar{a}_2\left[2r-\frac{1}{r}+r^3\right]\\&+&C_3\bar{a}_3\left[3r-\frac{2}{5r}+\frac{3r^3}{2}+r^5\right].
 	\end{eqnarray*}
\end{lemma}
\begin{proof}
The term  $z^n$ in $T_{f}T_{z+\overline{g}}(z^n)=T_{z+\overline{g}}T_{f}(z^n)$ appears both sides only from the expressions $$\left(T_{e^{-i\theta}f_{-1}}T_z+T_{e^{i\theta}f_1}T_{\bar{a}_1\bar{z}}+T_{e^{2i\theta}f_2}T_{\bar{a}_2\bar{z}^2}+T_{e^{3i\theta}f_3}T_{\bar{a}_3\bar{z}^3}\right)(z^n)$$ and $$\left(T_zT_{e^{-i\theta}f_{-1}}+T_{\bar{a}_1\bar{z}}T_{e^{i\theta}f_1}+T_{\bar{a}_2\bar{z}^2}T_{e^{2i\theta}f_2}+T_{\bar{a}_3\bar{z}^3}T_{e^{3i\theta}f_3}\right)(z^n),$$ respectively. So both sides are equal,  using the results of the previous lemmas and evaluating each term on both sides yields:
\begin{eqnarray}\label{ennik1}
(2n+2)\widehat{f_{-1}}(2n+3)-2n\widehat{f_{-1}}(2n+1)&=&C_1\bar{a}_1\left(\frac{2n+2}{2n+4}-\frac{2n}{2n+2}\right)\nonumber\\
&+&C_3\bar{a}_1^2\Bigg[3\left(\frac{2n+2}{2n+4}-\frac{2n}{2n+2}\right)\nonumber\\
&-&2\left(\frac{2n+2}{(2n+4)^2}-\frac{2n}{(2n+2)^2}\right)+\frac{2n+2}{2n+6}\nonumber\\&-&\frac{2n}{2n+4} \Bigg]\nonumber\\
&+&C_2\bar{a}_2\left(\frac{2n+2}{2n+6}-\frac{2n-2}{2n+2}\right)\nonumber\\&+&C_3\bar{a}_3\left(\frac{2n+2}{2n+8}-\frac{2n-4}{2n+2}\right).
\end{eqnarray}
We complexify the above equation by considering $z=2n-4$ and we use Remark \ref{periodic} to obtain
\begin{eqnarray}\label{ennik2}
(z+6)\widehat{f_{-1}}(z+7)-(z+4)\widehat{f_{-1}}(z+5)&=&C_1\bar{a}_1\left(\frac{z+6}{z+8}-\frac{z+4}{z+6}\right)\nonumber\\
&+&C_3\bar{a}_1^2\Bigg[3\left(\frac{z+6}{z+8}-\frac{z+4}{z+6}\right)\nonumber\\
&-&2\left(\frac{z+6}{(z+8)^2}-\frac{z+4}{(z+6)^2}\right)+\frac{z+6}{z+10}\nonumber\\&-&\frac{z+4}{z+8} \Bigg]
+C_2\bar{a}_2\left(\frac{z+6}{z+10}-\frac{z+2}{z+6}\right)\nonumber\\
&+&C_3\bar{a}_3\left(\frac{z+6}{z+12}-\frac{z}{z+6}\right).
\end{eqnarray}
We let $$F(z)=(z+4)\widehat{f_{-1}}(z+5)$$ and 
\begin{eqnarray*}G(z)&=&C_1\bar{a}_1\frac{z+4}{z+6}+C_3\bar{a}_1^2\left[3\frac{z+4}{z+6}-2\frac{z+4}{(z+6)^2}+\frac{z+4}{z+8} \right]+C_2\bar{a}_2\sum_{i=0}^1\frac{z+2i+2}{z+2i+6}\\
	&+&C_3\bar{a}_3\sum_{i=0}^2\frac{z+2i}{z+2i+6}.\end{eqnarray*}
Then, equation \eqref{ennik2} can be written as $F(z+2)-F(z)=G(z+2)-G(z)$. Thus, by Remark \ref{periodic}, there exists constant $C_{-1}$ such that $F(z)=C_{-1}+G(z)$. Applying partial fraction decomposition to the terms of  $G$ and using Remark \ref{periodic} implies that
\begin{eqnarray*}\widehat{r^5f_{-1}}(z)&=&C_{-1}\widehat{r^4}(z)+C_1\bar{a}_1\widehat{r^6}(z)+C_3\bar{a}_1^2\left[3\widehat{r^6}(z)+2\widehat{r^6\ln r}(z)+\widehat{r^8}(z)\right]\\ &+& C_2\bar{a}_2\left[2\widehat{r^6}(z)
	-\widehat{r^4}(z)+\widehat{r^8}(z) \right] +C_3\bar{a}_3\Big[3\widehat{r^6}(z)-2\widehat{r^4}(z)-\frac{1}{2}\widehat{r^4}(z)\\
	&+&\frac{3}{2}\widehat{r^8}(z)+\widehat{r^{10}}(z) \Big].\end{eqnarray*}
Hence,  \begin{eqnarray*}f_{-1}(r)&=&\frac{C_{-1}}{r}+C_1\bar{a}_1r+C_3\bar{a}_1^2\left[3r+2r\ln r+r^3\right]+C_2\bar{a}_2\left[2r-\frac{1}{r}+r^3\right]\\
	&+&C_3\bar{a}_3\left[3r-\frac{2}{5r}+\frac{3r^3}{2}+r^5\right].\end{eqnarray*}

\end{proof}
The main purpose of the following lemma is to evaluate the radial function $f_{-2}$. However, we will omit some of the lengthy calculations, as they are similar to those in the previous lemmas.
\begin{lemma}
\label{f-2}
Under the hypothesis of Theorem \ref{main}, we have that
\begin{eqnarray*}
f_{-2}(r)&=&\frac{C_{-2}}{r^2}-C_2\bar{a}_1^2\left(\frac{1}{r^2}-r^2\right)-C_3\bar{a}_1\bar{a}_2\Bigg(\frac{31}{4r^2}-
6r^2-2r^{4}-2r^2\ln r+\frac{1}{4r^6}\\
&+&\frac{1}{2r^4}-\frac{\ln r}{r^2}-\frac{1}{2} \Bigg)
+C_1\bar{a}_2r^2-C_2\bar{a}_3\left( \frac{3}{2r^2}-\frac{3}{2}r^2+\frac{1}{r^2}-r^4\right)\\&-&C_3\bar{a}_4\left( \frac{13}{3r^2}-2r^2-\frac{4}{3}r^4-r^6\right).
\end{eqnarray*}
\end{lemma}
\begin{proof}
The term $z^{n-1}$ appears in the left-hand side and right-hand side of equation \eqref{commute}  only from
$$\left(T_{e^{-2i\theta}f_{-2}}T_z+T_{f_0}T_{\bar{a}_1\bar{z}}+T_{e^{i\theta}f_1}T_{\bar{a}_2\bar{z}^2}+T_{e^{2i\theta}f_2}T_{\bar{a}_3\bar{z}^3}+T_{e^{3i\theta}f_3}T_{\bar{a}_4\bar{z}^4}\right)(z^n)$$ and 
$$\left(T_zT_{e^{-2i\theta}f_{-2}}+T_{\bar{a}_1\bar{z}}T_{f_0}+T_{\bar{a}_2\bar{z}^2}T_{e^{i\theta}f_1}+T_{\bar{a}_3\bar{z}^3}T_{e^{2i\theta}f_2}+T_{\bar{a}_4\bar{z}^4}T_{e^{3i\theta}f_3}\right)(z^n),$$
respectively. So both sides must be equal. Next, using Lemma \ref{mellin}, we evaluate each term on both sides and we obtain:
\begin{enumerate}
\item $T_{e^{-2i\theta}f_{-2}}T_z(z^n)=2n\widehat{f_{-2}}(2n+2)z^{n-1},$\\
\item \begin{eqnarray*}
T_{f_0}T_{\bar{a}_1\bar{z}}(z^n)&=&\bar{a}_1 \frac{(2n)^2}{2n+2}\widehat{f_0}(2n)z^{n-1}\\
&=& \Bigg[C_0\bar{a}_1\frac{2n}{2n+2}+C_2\bar{a}_1^2\Big( \frac{2n}{2n+2}-\frac{2}{2n+2}+\frac{(2n)^2}{(2n+2)^2}\Big)\\
&+&C_3\bar{a}_1\bar{a}_2\left(-\frac{4}{2n+2}+2\frac{(2n)^2}{(2n+2)^2}+\frac{(2n)^2}{(2n+2)(2n+4)} \right) \Bigg]z^{n-1},\\
\end{eqnarray*}
\item \begin{eqnarray*}
 T_{e^{i\theta}f_1}T_{\bar{a}_2\bar{z}^2}(z^n)&=&\bar{a}_2 \frac{2n(2n-2)}{2n+2}\widehat{f_1}(2n-1)z^{n-1}\\
&=& \Bigg[C_1\bar{a}_2\frac{2n-2}{2n+2}+C_3\bar{a}_1\bar{a}_2\Big( \frac{2n(2n-2)}{(2n+2)^2}+3\frac{2n-2}{2n+2}-2\frac{2n-2}{2n(2n+2)}\\
&-&\frac{2n}{2n+2}\Big) \Bigg]z^{n-1},\\
\end{eqnarray*}
\item $\displaystyle T_{e^{2i\theta}f_2}T_{\bar{a}_3\bar{z}^3}(z^n)=C_2\bar{a}_3\frac{2n-4}{2n+2}z^{n-1}$,\\
\item $\displaystyle T_{e^{3i\theta}f_3}T_{\bar{a}_4\bar{z}^4}(z^n)=C_3\bar{a}_4\frac{2n-6}{2n+2}z^{n-1}$,\\
\item $T_zT_{e^{-2i\theta}f_{-2}}(z^n)=(2n-2)\widehat{f_{-2}}(2n)z^{n-1}$,\\
\item \begin{eqnarray*}
T_{\bar{a}_1\bar{z}}T_{f_0}(z^n)&=&2n\bar{a}_1(2n+2)\widehat{f_0}(2n+2)\widehat{r}(2n+1)z^{n-1}\\
&=& \Bigg[C_0\bar{a}_1\frac{2n}{2n+2}+C_2\bar{a}_1^2\left( \frac{2n}{2n+2}-\frac{4n}{(2n+2)^2}+\frac{2n}{2n+4}\right)\\
&+&C_3\bar{a}_1\bar{a}_2\left(-\frac{8n}{(2n+2)^2}+\frac{4n}{2n+4}+\frac{2n}{2n+6} \right) \Bigg]z^{n-1},\\
\end{eqnarray*}
\item \begin{eqnarray*} T_{\bar{a}_2\bar{z}^2}T_{e^{i\theta}f_1}(z^n)&=& \Bigg[C_1\bar{a}_2\frac{2n}{2n+4}+C_3\bar{a}_1\bar{a}_2\Big( \frac{2n}{2n+6}+\frac{6n}{2n+4}-\frac{4n}{(2n+4)^2}\\
&-&\frac{2n}{2n+2}\Big) \Bigg]z^{n-1},
\end{eqnarray*}
    \item $\displaystyle T_{\bar{a}_3\bar{z}^3}T_{e^{2i\theta}f_2}(z^n)=C_2\bar{a}_3\frac{2n}{2n+6}z^{n-1}$,\\
    \item $\displaystyle T_{\bar{a}_4\bar{z}^4}T_{e^{3i\theta}f_3}(z^n)=C_3\bar{a}_4\frac{2n}{2n+8}z^{n-1}.$
\end{enumerate}
We substitutes the ten quantities above into both sides, equate them, complexify the expression by setting $z=2n-6$, and then use Remark \ref{periodic} to obtain the following:
\begin{eqnarray}
\label{complexified}
(z+6)\widehat{f_{-2}}(z+8)-(z+4)\widehat{f_{-2}}(z+6)&=&C_2\bar{a}_1^2\left( \frac{z+6}{z+8}-\frac{2(z+6)}{(z+8)^2}+\frac{z+6}{z+10}\right)\nonumber\\
&-&C_2\bar{a}_1^2\left( \frac{z+6}{z+8}-\frac{2}{z+8}+\frac{(z+6)^2}{(z+10)^2}\right)\nonumber\\
&+&C_3\bar{a}_1\bar{a}_2\left( \frac{-4(z+6)}{(z+8)^2}+\frac{2(z+6)}{z+10}+\frac{z+6}{z+12}\right)\nonumber\\
&-&C_3\bar{a}_1\bar{a}_2\left( \frac{-4}{z+8}+\frac{2(z+6)^2}{(z+8)^2}+\frac{(z+6)^2}{(z+8)(z+10)}\right)\nonumber\\
&+&C_1\bar{a}_2\frac{z+6}{z+10}-C_1\bar{a}_2\frac{z+4}{z+8}\nonumber\\
&+&C_3\bar{a}_1\bar{a}_2\left( \frac{z+6}{z+12}+\frac{3(z+6)}{z+10}-\frac{2(z+6)}{(z+10)^2}\right)\nonumber\\
&-&C_3\bar{a}_1\bar{a}_2\left( \frac{(z+4)(z+6)}{(z+8)^2}+\frac{3(z+4)}{z+8}-\frac{2(z+4)}{z(z+8)}\right)\nonumber\\
&+&C_2\bar{a}_3\frac{z+6}{z+12}-C_2\bar{a}_3\frac{z+2}{z+8}\nonumber\\
&+& C_3\bar{a}_4\frac{z+6}{z+8}-C_3\bar{a}_4\frac{z}{z+8}.\\
\nonumber\end{eqnarray}
Now, let $F(z)=(z+4)\widehat{f_{-2}}(z+6)$ and define $G(z)$ as $\displaystyle G(z)=\sum_{i=1}^9G_i(z)$, where $G_i(z)$ are given by
\begin{enumerate}
\item  $\displaystyle G_1(z)=-4C_2\bar{a}_1^2\frac{1}{z+8}$,\\

\item $\displaystyle G_2(z)=-6C_3\bar{a}_1\bar{a}_2\sum_{i=0}^1\frac{1}{z+2i+8}$,\\

\item $\displaystyle G_3(z)=C_1\bar{a}_2\frac{z+4}{z+8}$,\\

\item $\displaystyle G_4(z)=-6C_3\bar{a}_1\bar{a}_2\sum_{i=0}^1\frac{1}{z+2i+8}$,\\

\item $\displaystyle G_5(z)=-14C_3\bar{a}_1\bar{a}_2\frac{1}{z+8}$,\\

\item $\displaystyle G_6(z)=8C_3\bar{a}_1\bar{a}_2\frac{1}{(z+8)^2}$,\\

\item $\displaystyle G_7(z)=-C_3\bar{a}_1\bar{a}_2\sum_{i=0}^3\frac{1}{z+2i}$,\\

\item $\displaystyle G_8(z)=-6C_2\bar{a}_3\sum_{i=0}^1\frac{1}{z+2i+8}$,\\

\item $\displaystyle G_9(z)=-6C_3\bar{a}_4\sum_{i=0}^2\frac{1}{z+2i+8}$.\\
\end{enumerate}
Thus, equation \eqref{complexified} simplifies to $F(z+2)-F(z)=G(z+2)-G(z)$. Therefore, by Remark \ref{periodic}, there exists a constant $C_{-2}$ such that $F(z)=C_{-2}+G(z)$, which is equivalent to $$(z+4)\widehat{f_{-2}}(z+6)=C_{-2}+G(z).$$ Now, dividing both sides by $(z+4)$ and expanding $G(z)$ using sums of partial fractions, we obtain:
\begin{eqnarray*}
\label{mnayykah}
\widehat{f_{-2}}(z+6)&=&\frac{C_{-2}}{z+4}-C_{2}\bar{a}_1^2\left( \frac{1}{z+4}-\frac{1}{z+8}\right)-\frac{3}{2}C_{3}\bar{a}_1\bar{a}_2\left(\frac{1}{z+4}-\frac{1}{z+8}\right)\\
&-&C_{3}\bar{a}_1\bar{a}_2\left( \frac{1}{z+4}-\frac{1}{z+10}\right)
+C_1\bar{a}_2\frac{1}{z+8}-\frac{3}{2}C_{3}\bar{a}_1\bar{a}_2\left( \frac{1}{z+4}-\frac{1}{z+8}\right)\\
&-&C_{3}\bar{a}_1\bar{a}_2\left( \frac{1}{z+4}-\frac{1}{z+10}\right)-\frac{7}{2}C_{3}\bar{a}_1\bar{a}_2\left( \frac{1}{z+4}-\frac{1}{z+8}\right)\nonumber\\
&+&\frac{1}{2}C_3\bar{a}_1\bar{a}_2 \frac{1}{z+4}-\frac{1}{2}C_3\bar{a}_1\bar{a}_2 \frac{1}{z+8}-2C_3\bar{a}_1\bar{a}_2 \frac{1}{(z+8)^2}\\
&-&\frac{1}{4}C_{3}\bar{a}_1\bar{a}_2\left( \frac{1}{z}-\frac{1}{z+4}\right)-\frac{1}{2}C_{3}\bar{a}_1\bar{a}_2\left( \frac{1}{z+2}-\frac{1}{z+4}\right)\nonumber\\
&-&C_{3}\bar{a}_1\bar{a}_2 \frac{1}{(z+4)^2}-\frac{1}{2}C_{3}\bar{a}_1\bar{a}_2\left( \frac{1}{z+4}-\frac{1}{z+6}\right)-\frac{3}{2}C_{2}\bar{a}_3\left( \frac{1}{z+4}-\frac{1}{z+8}\right)\\
&-&C_{2}\bar{a}_3\left( \frac{1}{z+4}-\frac{1}{z+10}\right)
-C_{3}\bar{a}_4\left( \frac{1}{z+4}-\frac{1}{z+8}\right)\\
&-&\frac{4}{3}C_{3}\bar{a}_4\left( \frac{1}{z+4}-\frac{1}{z+10}\right)-C_{3}\bar{a}_4\left( \frac{1}{z+4}-\frac{1}{z+12}\right).
\end{eqnarray*}
Finally, applying Remark \ref{periodic}, the equation above implies
\begin{eqnarray*}
f_{-2}(r)&=&\frac{C_{-2}}{r^2}-C_2\bar{a}_1^2\left(\frac{1}{r^2}-r^2\right)-C_3\bar{a}_1\bar{a}_2\Bigg(\frac{31}{4r^2}-
6r^2-2r^{4}-2r^2\ln r+\frac{1}{4r^6}\\
&+&\frac{1}{2r^4}-\frac{\ln r}{r^2}-\frac{1}{2} \Bigg)
+C_1\bar{a}_2r^2-C_2\bar{a}_3\left( \frac{3}{2r^2}-\frac{3}{2}r^2+\frac{1}{r^2}-r^4\right)\\
&-&C_3\bar{a}_4\left( \frac{13}{3 r^2}-2r^2-\frac{4}{3}r^4-r^6\right).
\end{eqnarray*}
\end{proof}
 We observe that $f_{-2}$, obtained in the previous lemma, belongs to $L^1([0,1),rdr)$ if and only if $C_{-2}=0$, $C_2=0$ and $C_3=0$. Thus, by Remark \ref{rmk1}, Lemma \ref{f0} and Lemma \ref{f-1}, we establish the following:
\begin{enumerate}
\item $T_{e^{3i\theta}f_3}=0$,\\
\item $T_{e^{2i\theta}f_2}=0$,\\
\item $T_{e^{i\theta}f_1}=C_1T_z$,\\
\item $T_{f_0}=C_0$,\\
\item $f_{-1}=\frac{C_{-1}}{r}+C_1\bar{a}_1r$,\\
\item $T_{e^{-2i\theta}f_{-2}}=C_1T_{\bar{a}_2\bar{z}^2}$.\\
\end{enumerate}
This implies that $N=1$ in the polar decomposition of the symbol $f$ in Theorem \ref{main}, and that $f(re^{i\theta})=\displaystyle{\sum_{k=-\infty}^{1}e^{ik\theta}f_{k}(r)}$.

\begin{lemma}
\label{f-1}
Under the hypothesis of Theorem \ref{main}, we have that  $$f_{-1}(r)=C_1\bar{a}_1r\ \textrm{ and }\  f_{-3}(r)=C_1\bar{a}_3r^3.$$
\end{lemma}
\begin{proof}
In equation \eqref{commute1}, the term $\bar{z}^{n+2}$ arises from:
\begin{eqnarray}
\label{bar2}
\left(T_{e^{-i\theta}f_{-1}}T_{\bar{a}_1\bar{z}}+T_{f_0}T_{\bar{a}_2\bar{z}^2}+T_{e^{i\theta}f_1}T_{\bar{a}_3\bar{z}^3}+T_{e^{-3i\theta}f_{-3}}T_{z}\right)(\bar{z}^n)&=&\nonumber\\
\left(T_{\bar{a}_1\bar{z}}T_{e^{-i\theta}f_{-1}}+T_{\bar{a}_2\bar{z}^2}T_{f_0}+T_{\bar{a}_3\bar{z}^3}T_{e^{i\theta}f_1}+T_{z}T_{e^{-3i\theta}f_{-3}}\right)(\bar{z}^n)&&
\end{eqnarray}
By using Lemma \ref{mellin}, we evaluate each term in equation \eqref{bar2} and we obtain
\begin{enumerate}
\item \begin{eqnarray*} T_{e^{-i\theta}f_{-1}}T_{\bar{a}_1\bar{z}}(\bar{z}^n)&=&\bar{a}_12(n+3)\widehat{f_{-1}}(2n+5)\bar{z}^{n+2}\\
	&=&\bar{a}_1(2n+6)\left[ \frac{C_{-1}}{2n+4}+\frac{C_1\bar{a}_1}{2n+6}\right]\bar{z}^{n+2},\end{eqnarray*}\\

    \item Since $f_0(r)=C_0$,  $T_{f_0}T_{\bar{a}_2\bar{z}^2}(\bar{z}^n)=T_{\bar{a}_2\bar{z}^2}T_{f_0}(\bar{z}^n)=C_0\bar{z}^{n+2}$,\\
    \item \begin{eqnarray*} T_{e^{i\theta}f_{1}}T_{\bar{a}_3\bar{z}^3}(\bar{z}^n)&=&C_1\bar{a}_32(n+3)\widehat{r}(2n+7)\bar{z}^{n+2}\\&=&C_1\bar{a}_3\frac{2n+6}{2n+8}\bar{z}^{n+2},\end{eqnarray*}

    \item $\displaystyle T_{e^{-3i\theta}f_{-3}}T_{z}(\bar{z}^n)=\frac{2n(2n+6)}{2n+2}\widehat{f_{-3}}(2n+3)\bar{z}^{n+2}$,\\

    \item \begin{eqnarray*} T_{\bar{a}_1\bar{z}}T_{e^{-i\theta}f_{-1}}(\bar{z}^n)&=&\bar{a}_1(2n+4)\widehat{f_{-1}}(2n+3)\bar{z}^{n+2}\\
    	&=&\bar{a}_1(2n+4)\left[ \frac{C_{-1}}{2n+2}+\frac{C_1\bar{a}_1}{2n+4}\right]\bar{z}^{n+2},\end{eqnarray*}\\

      \item $\displaystyle T_{\bar{a}_3\bar{z}^3}T_{e^{i\theta}f_{1}}(\bar{z}^n)=C_1\bar{a}_3\frac{2n}{2n+2}\bar{z}^{n+2}$,\\

      \item $\displaystyle T_{z}T_{e^{-3i\theta}f_{-3}}(\bar{z}^n)=(2n+6)\widehat{f_{-3}}(2n+5)\bar{z}^{n+2}$.\\

\end{enumerate}
Substituting the above quantities into equation \eqref{bar2} and rearranging them yields
\begin{eqnarray*}(2n+6)\widehat{f_{-3}}(2n+5)-\frac{2n(2n+6)}{2n+2}\widehat{f_{-3}}(2n+3)&=&C_{-1}\bar{a}_1\frac{2n+6}{2n+4}-C_{-1}\bar{a}_1\frac{2n+4}{2n+2}\\
	&+&C_{1}\bar{a}_3\frac{2n+6}{2n+8}-C_{1}\bar{a}_3\frac{2n}{2n+2},\end{eqnarray*}
which is equivalent to
\begin{eqnarray}\label{f_3}(2n+2)\widehat{f_{-3}}(2n+5)-2n\widehat{f_{-3}}(2n+3)&=&C_{-1}\bar{a}_1\frac{2n+2}{2n+4}-C_{-1}\bar{a}_1\frac{2n+4}{2n+6}\nonumber\\&+&C_{1}\bar{a}_3\frac{2n+2}{2n+8}-C_{1}\bar{a}_3\frac{2n}{2n+6}.\end{eqnarray}
We set $z=2n$ to complexify equation \eqref{f_3} and obtain:
$$(2z+2)\widehat{f_{-3}}(z+5)-z\widehat{f_{-3}}(z+3)=C_{-1}\bar{a}_1\frac{z+2}{z+4}-C_{-1}\bar{a}_1\frac{z+4}{z+6}+C_{1}\bar{a}_3\frac{z+2}{z+8}-C_{1}\bar{a}_3\frac{z}{z+6}.$$
We let $F(z)=z\widehat{f_{-3}}(z+3)$ and $G(z)$ be defined as $\displaystyle G(z)=C_{1}\bar{a}_3\frac{z}{z+6}-C_{-1}\bar{a}_1\frac{z+2}{z+4}$. Then the above equation simplifies to $$F(z+2)-F(z)=G(z+2)-G(z).$$
Therefore, by Remark \ref{periodic}, there exists a constant $C_{-3}$ such that $F(z)=C_{-3}+G(z)$. Hence, we have
\begin{eqnarray*}
\widehat{f_{-3}}(z+3)&=&\frac{C_{-3}}{z}+C_{1}\bar{a}_3\frac{1}{z+6}-C_{-1}\bar{a}_1\frac{z+2}{z(z+4)}\\
&=&\frac{C_{-3}}{z}+C_{1}\bar{a}_3\frac{1}{z+6}-C_{-1}\bar{a}_1\left( \frac{1}{2z}+\frac{1}{2(z+4)}\right).\\
\end{eqnarray*}
Using Remark \ref{periodic}, we deduce that
\begin{equation*}
\widehat{r^3f_{-3}}(z)=C_{-3}\hat{1}(z)+C_{1}\bar{a}_3\widehat{r^6}(z)-C_{-1}\bar{a}_1\left( \frac{1}{2}\hat{1}(z)+\frac{1}{2}\widehat{r^4}(z)\right).\\
\end{equation*}
Therefore,
$$f_{-3}(r)=\frac{C_{-3}}{r^3}+C_{1}\bar{a}_3r^3- \frac{1}{2}C_{-1}\bar{a}_1\left( \frac{1}{r^3}+r\right).$$
Clearly, $f_{-3}(r)$ belongs to $L^1([0,1),rdr)$ if and only if $C_{-3}=0$ and $C_{-1}=0$. In this case, we have $f_{-1}(r)=C_1\bar{a}_1r$ and $f_{-3}(r)=C_1\bar{a}_3r^3$.
\end{proof}

We now proceed with the computation of $f_{-4}$.
\begin{lemma}
\label{f-4}
Under the hypothesis of Theorem \ref{main}, we have that $f_{-4}(r)=C_1\bar{a}_4r^4.$
\end{lemma}
\begin{proof}
In equation \eqref{commute1}, the term $\bar{z}^{n+3}$ arises from the following expression:
\begin{eqnarray}
\label{bar1}
\left(T_{e^{-4i\theta}f_{-4}}T_{z}+T_{f_0}T_{\bar{a}_3\bar{z}^3}+T_{e^{i\theta}f_1}T_{\bar{a}_4\bar{z}^4}+T_{e^{-i\theta}f_{-1}}T_{\bar{a}_2\bar{z}^2}\right)(\bar{z}^n)&=&\nonumber\\
\left(T_{z}T_{e^{-4i\theta}f_{-4}}+T_{\bar{a}_3\bar{z}^3}T_{f_0}+T_{\bar{a}_4\bar{z}^4}T_{e^{i\theta}f_1}+T_{\bar{a}_2\bar{z}^2}T_{e^{-i\theta}f_{-1}}\right)(\bar{z}^n)&&
\end{eqnarray}
We use Lemma \ref{mellin} to evaluate each term appearing in equation \eqref{bar1}, and we obtain:
\begin{enumerate}
\item $\displaystyle T_{e^{-4i\theta}f_{-4}}T_{z}(\bar{z}^n)=\frac{2n(2n+8)}{2n+2}\widehat{f_{-4}}(2n+4)\bar{z}^{n+3},$\\
    \item since $f_0(r)=C_0$,  $T_{f_0}T_{\bar{a}_3\bar{z}^3}(\bar{z}^n)=T_{\bar{a}_3\bar{z}^3}T_{f_0}(\bar{z}^n)=C_0\bar{z}^{n+3}$,\\
    \item $\displaystyle T_{e^{i\theta}f_{1}}T_{\bar{a}_4\bar{z}^4}(\bar{z}^n)=C_1\bar{a}_42(n+4)\widehat{r}(2n+9)\bar{z}^{n+3}=C_1\bar{a}_4\frac{2n+8}{2n+10}\bar{z}^{n+3}$,\\
    \item $\displaystyle T_{e^{-i\theta}f_{-1}}T_{\bar{a}_2\bar{z}^2}(\bar{z}^n)=C_1\bar{a}_1\bar{a}_2\bar{z}^{n+3},$\\
     \item $\displaystyle T_{z}T_{e^{-4i\theta}f_{-4}}(\bar{z}^n)=(2n+8)\widehat{f_{-4}}(2n+6)\bar{z}^{n+3}$,\\
      \item $\displaystyle T_{\bar{a}_4\bar{z}^4}T_{e^{i\theta}f_{1}}(\bar{z}^n)=C_1\bar{a}_4\frac{2n}{2n+2}\bar{z}^{n+3}$.\\
\end{enumerate}
After substituting the above terms in equation \eqref{bar1} and rearranging them, we obtain
$$(2n+2)\widehat{f_{-4}}(2n+6)-2n\widehat{f_{-4}}(2n+4)=C_1\bar{a}_4\frac{2n+2}{2n+10}-C_1\bar{a}_4\frac{2n}{2n+8}.$$
Thus, by setting $z=2n$, the above equation becomes
$$F(z+2)-F(z)=G(z+2)-G(z),$$
where $F(z)=z\widehat{f_{-4}}(z+4)$ and $G(z)=C_1\bar{a}_4\frac{z}{z+8}$. Remark \ref{periodic} implies the existence of a constant $C_{-4}$ such that $F(z)=C_{-4}+G(z)$. Hence,
$$\widehat{r^4f_{-4}}=C_{-4}\hat{1}(z)+C_1\bar{a}_4\widehat{r^8}(z).$$
Therefore, we deduce that $\displaystyle f_{-4}(r)=\frac{C_{-4}}{r^4}+C_1\bar{a}_4r^4$. Clearly $f_4$ belongs to  $L^1([0,1),rdr)$ if and only if $C_{-4}=0$. Finally, we must have $f_{-4}(r)=C_1\bar{a}_4r^4$.
\end{proof}
Using the same technique as in the previous lemmas, we establish the following by induction.
\begin{lemma}
If equations \ref{commute} and \ref{commute1} are satisfied, then for all $k\geq 1$, we have $f_{-k}(r)=C_1\bar{a}_kr^k$ 
\end{lemma}
\begin{proof}
By Lemma \ref{f-1}, we have $f_{-1}(r)=C_1\bar{a}_1r$, which establish the base case. Now, assume that the formula holds for some $k\geq 1$, i.e., $f_{-k}(r)=C_1\bar{a}_kr^k$. Following a similar argument as in the proof of Lemma  \ref{f-4}, we obtain $\displaystyle f_{-(k+1)}(r)=\frac{C_{-(k+1)}}{r^{k+1}}+C_1\bar{a}_{k+1}r^{k+1}$. For $f_{-(k+1)}$ to belong to $L^1([0,1),rdr)$, it must satisfy $C_{-(k+1)}=0$. Thus, we conclude that $f_{-(k+1)}(r)=C_1\bar{a}_{k+1}r^{k+1}$, which completes the induction.
\end{proof}

\section{Conclusion}
Combining all the results from the previous section, we conclude that the symbol $f$ can now be expressed as $$f(z)=C_1z+C_0+C_1\sum_{n=1}^\infty \bar{a_n}\bar{z}^n.$$ This, in turn, implies that the Toeplitz operator $T_f$ takes the form $$T_f=C_1T_u+C_0I,$$ where $I$ denotes the identity operator. This completes the proof of Theorem \ref{main}.\\

\noindent{\bf{Final Remark:}}   
The results in this paper describe bounded Toeplitz operators, with truncated symbols, that commute with \( T_{z+\bar{g}} \), where \( g \) is an analytic function. It is worth noting that the analytic part \( z \) in the symbol \( z+\bar{g} \) can be replaced by \( z^n \) or even a polynomial in \( z \), and the same proof techniques can still be applied. However, this generalization comes at the cost of significantly more involved calculations, which can become tedious and lengthy. 

{\small


\begin{thebibliography}{999}
	
	    \bibitem{chen} Chen, Y., Koo, H. \& Lee, Y.J. \textit{Ranks of commutators of Toeplitz operators on the harmonic Bergman space}. Integr. Equ. Oper. Theory {\bf{75}}, 31–38 (2013). https://doi.org/10.1007/s00020-012-2020-6
	    \bibitem{choe} Choe, Boo Rim and Lee, Young Joo. \textit{Commuting Toeplitz operators on the harmonic Bergman space}. Michigan Math. J. {\bf{49}} (1999), no. 1, 163-174.
	     	    \bibitem{choe1} Choe, B.R., Lee, Y.J. \textit{Commutants of Analytic Toeplitz Operators on the Harmonic Bergman Space}. Integr. equ. oper. theory {\bf{50}}, 559–564 (2004). https://doi.org/10.1007/s00020-004-1338-0
	    \bibitem{choe2} B. R. Choe, Y. J. Lee and K. Na. \textit{Toeplitz operators on harmonic Bergman spaces}. Nagoya Math. J. Vol. {\bf{174}} (2004), 165-186.
	    \bibitem{ding} Ding, Xuanhao.
	    \textit{A question of Toeplitz operators on the harmonic Bergman space}.
	    J. Math. Anal. Appl. {\bf{344}} (2008), no. 1, 367-372.
	    \bibitem{ding2} Ding, Qian; Hu, Yinyin; Liu, Liu; Lu, Yufeng.
	    \textit{Zero products and finite rank of Toeplitz operators on the harmonic Bergman space}. J. Math. Res. Appl. {\bf{37}} (2017), no. 3, 325-334.
	    \bibitem{dong1} Dong, Xing-Tang; Liu, Congwen; Zhou, Ze-Hua.
	    \textit{Quasihomogeneous Toeplitz operators with integrable symbols on the harmonic Bergman space}. Bull. Aust. Math. Soc. {\bf{90}} (2014), no. 3, 494-503.
	\bibitem{dong2} X. T. Dong and Z. H. Zhou, \textit{Products of Toeplitz operators on the harmonic Bergman space}, Proc. Amer. Math. Soc. {\bf{138}} (2010), 1765-1773.
	\bibitem{dong3} X. T. Dong and Z. H. Zhou, \textit{Commuting quasihomogeneous Toeplitz operators on the harmonic
		Bergman space}, Complex Anal. Oper. Theory {\bf{7}} (2013), 1267-1285.
	\bibitem{dong4} X. T. Dong and Z. H. Zhou, \textit{Product equivalence of quasihomogeneous Toeplitz operators on the
		harmonic Bergman space}, Studia Math. {\bf{219}} (2013), 163-175.
		\bibitem{guan} Guan, H.Y., Lu, Y.F. \textit{Algebraic properties of Toeplitz and small Hankel operators on the harmonic Bergman space}. Acta. Math. Sin.-English Ser. {\bf{30}}, 1395-1406 (2014). https://doi.org/10.1007/s10114-014-3276-3
		\bibitem{guo} Guo, Kunyu; Zheng, Dechao.
		\textit{Toeplitz algebra and Hankel algebra on the harmonic Bergman space}.
		J. Math. Anal. Appl. {\bf{276}} (2002), no. 1, 213-230.
			
		\bibitem{kong} Kong, Linghui; Qu, Shuang; Tong, Shan.
		 \textit{Commuting and semi-commuting Toeplitz operators on the weighted harmonic Bergman space}.
		Oper. Matrices {\bf{15}} (2021), no. 1, 163-174.
		\bibitem{le} Le, Trieu. \textit{
		Toeplitz operators on radially weighted harmonic Bergman spaces}.
		J. Math. Anal. Appl. {\bf{396}} (2012), no. 1, 164-172.	
		\bibitem{le2} Le, Trieu and Tikaradze, Akaki. \textit{Commutants of Toeplitz operators with harmonic symbols}. New York J. Math. {\bf{23}} (2017), 1723-1731.	
		\bibitem{lz} Louhichi, I., Zakariasy, L. \textit{Quasihomogeneous Toeplitz operators on the harmonic Bergman space}. Arch. Math. (Basel) {\bf{98}}, 49-60 (2012). https://doi.org/10.1007/s00013-011-0346-y
		\bibitem{lzf} Louhichi, Issam, Randriamahaleo, Fanilo and Zakariasy, Lova. \textit{On the Commutativity of a Certain Class of Toeplitz Operators}. Concr. Oper. {\bf{2}} (2015), no. 1, 1-7. https://doi.org/10.2478/conop-2014-0001
		\bibitem{Rem} R. Remmert, \textit{Classical Topics in complex function theory}, Graduate Texts in Mathematics \textbf{172}, Springer, New York, 1998.
	
	



		
	\end{thebibliography}
\end{document}